\documentclass[11pt]{article}
\usepackage[T1]{fontenc}
\usepackage{lmodern,amsmath,amssymb,amsthm,tikz-cd}
\usepackage{tikz}
\usetikzlibrary{arrows.meta,decorations.pathreplacing}
\usepackage[a4paper,margin=26mm]{geometry}
\usepackage[hidelinks]{hyperref}
\newcommand{\g}{\mathfrak g}
\newcommand{\db}{\mathfrak d}
\newcommand{\TT}{\mathbb T}
\newcommand{\HH}{\mathcal H}
\newcommand{\pr}{\operatorname{pr}}
\newcommand{\Gr}{\operatorname{Gr}}
\newcommand{\Lie}{\operatorname{Lie}}
\newcommand{\ad}{\operatorname{ad}}
\newcommand{\fus}{\mathrm{fus}}
\newcommand{\src}{\mathsf{s}}
\newcommand{\tgt}{\mathsf{t}}
\newcommand{\sect}[1]{\S\ref{#1}}

\makeatletter
\let\l@paragraph\l@subsection
\def\toclevel@paragraph{2}
\renewcommand{\@seccntformat}[1]{\S\csname the#1\endcsname\quad}
\makeatother
\newtheorem{thm}{Theorem}[section]
\newtheorem{prop}{Proposition}[section]
\newtheorem{cor}[prop]{Corollary}
\numberwithin{equation}{section}
\theoremstyle{definition}
\newtheorem{exa}[thm]{Example}
\theoremstyle{remark}
\newtheorem{rema}[thm]{Remark}
\title{Lagrangian structures, \(r\)-matrices and Poisson brackets}
\author{Daniel Álvarez\thanks{
Instituto de Matemática Pura e Aplicada (IMPA),
Rio de Janeiro, Brazil.
Email: \href{mailto:uerbum@impa.br}{\texttt{uerbum@impa.br}}.}}
\date{}
\begin{document}
\maketitle
\begin{abstract}
    We give a global conceptual viewpoint of two constructions of Poisson structures involving $r$-matrices: the first one concerns coupling of generalized dynamical $r$-matrices with actions of quasi-triangular Lie bialgebras, the second one is given by the mixed products of Poisson spaces, which combine two Poisson actions of a quasi-triangular Lie bialgebra to create a new one. When the input Lie algebroids are integrable, we construct symplectic groupoids for both types of Poisson structures by taking fibred products of Lie groupoid morphisms equipped with shifted lagrangian structures. We illustrate our results with explicit examples coming from Lie theory and moduli spaces of flat bundles over surfaces.   
\end{abstract}
\tableofcontents

\section{Introduction}\label{intro}

Moment maps are central in symplectic geometry, given their numerous applications to reduction, dynamics, convexity and equivariant cohomology. In Poisson geometry, the very concept of moment map has been extended in order to accommodate various remarkable examples of possibly Lie group or homogeneous space valued moment maps \cite{AKS,AMM,LuMomentum}. Given the non linear nature of these generalized moment maps, their existence and equivariance may not be treated by simply using Lie algebra cohomology as in the case of Lie group actions by symplectomorphisms. So it is remarkable that the Lie theoretic approach to Poisson geometry provides a unifying theory for a very wide family of such non linear moment maps.

In fact, Poisson manifolds may be regarded as infinitesimal objects, in the sense that their cotangent bundles are naturally equipped with a Lie algebroid structure. The underlying Lie bracket is determined by the Poisson bracket \(\{\cdot,\cdot\}\) via the formula
\[ [df,dg]= d\{f,g\}, \]
for all pairs of functions $f,g$ on the manifold; in fact, we may extend this expression to all 1-forms by using the Leibniz rule. The global objects corresponding to Poisson manifolds are known as symplectic groupoids, these are symplectic manifolds possessing a compatible Lie groupoid structure, encoding the Lie bracket described above. Because of this reason, symplectic groupoids are then termed integrations of Poisson manifolds. It is by looking at these objects that we obtain a very general theory of equivariant moment maps suitable for performing moment map reduction in most important situations.

As a first instance of this phenomenon, we have that Poisson manifolds equipped with Poisson actions of Poisson Lie groups
give rise, upon integration, to Hamiltonian actions. For example, let
\(G\) be a connected complete Poisson Lie group and \(G^*\) its simply
connected dual group. A Poisson action of \(G\) on an integrable Poisson
manifold \((M,\pi_M)\) lifts to a Hamiltonian action on its source-simply
connected symplectic groupoid \(\Sigma(M)\), with moment map
\[J:\Sigma(M)\to G^*\] 
\cite[Thm.~1]{FIP}. In particular, for an action by Poisson
diffeomorphisms we obtain a classical moment map with values in the dual of a Lie
algebra \cite[Thm.~1.1]{FOR}.

This picture reflects the fact that a Poisson action of a Lie bialgebra
\((\g,\g^*)\) on \((M,\pi_M)\) may be encoded in a Lie algebroid morphism
\[T^*_{\pi_M}M\to\g^*\] 
given by the transpose of the infinitesimal action in the convention of
\cite[Prop.~1.8]{FIP}. In the convention of \cite[Ex.~3.19]{PartI}
used below, this morphism is \(-\rho^*\), where \(\rho^*\) is the transpose
of the infinitesimal action \(\rho:M\times\g\to TM\). Here \(T^*_{\pi_M}M\) denotes the cotangent Lie
algebroid, with anchor \( \alpha \mapsto \pi_M^\sharp\alpha\). We may
additionally encode the $\mathfrak{g}$-action on $M$ in a Lie algebroid structure on
\(A_M=T^*M\oplus(M\times\g)\), together with a Lie algebroid morphism to the
Drinfeld double \(\db=\g\bowtie\g^*\) of the Lie bialgebra \cite{Lu}.
The anchor \(\mathsf a:A_M\to TM\) and the morphism
\(\varphi:A_M\to\db\) are given by
\[
 \mathsf a(\alpha,u)=\pi_M^\sharp\alpha+\rho(u),\qquad
 \varphi(\alpha,u)=(u,-\rho^*\alpha).
\]
From this viewpoint, the global counterpart of \(\varphi\), when
\(A_M\) is integrable, is a {\em multiplicative \(D\)-valued moment map}
\[
                 \Phi:\HH\longrightarrow D ,
\]
where \(\HH\rightrightarrows M\) is the source-simply connected integration of \(A_M\), \(D\) integrates
\(\db\), and $\Phi$ is a Lie groupoid morphism which lifts \(\varphi\) and is \(\g\times\g\)-equivariant \cite[Thm.~5.3]{PartI}. The groupoid \(\HH\) carries a
two-form \(\beta\) that fails to be a symplectic form compatible with $\HH$ in a way controlled by the moment map $\Phi$. Indeed, letting \(\Omega_D\) and
\(\Theta_D\) be, respectively, the Polyakov--Wiegmann two-form on \(D\times D\)
and the Cartan three-form on \(D\), \(\beta\) satisfies
\[
 \ker \beta \cap \ker T\Phi =0,\qquad \pr_1^*\beta+\pr_2^*\beta-m_{\HH}^*\beta=(\Phi\times\Phi)^*\Omega_D,
 \qquad d\beta=-\Phi^*\Theta_D ,
\]
where \(m_{\HH}\) is the groupoid multiplication and
\(\pr_i\) are the projections from the space of composable pairs. In the modern language of shifted symplectic geometry, $\beta$ is an example of a {\em (2-shifted) lagrangian structure on $\Phi$} \cite{PTVV}. The main utility of such $D$-valued moment maps (and of lagrangian structures in general) is in performing a generalized version of symplectic reduction involving fibred products and quotients by Lie group actions, which accounts for the construction of many symplectic groupoids of interest \cite{PartI}.

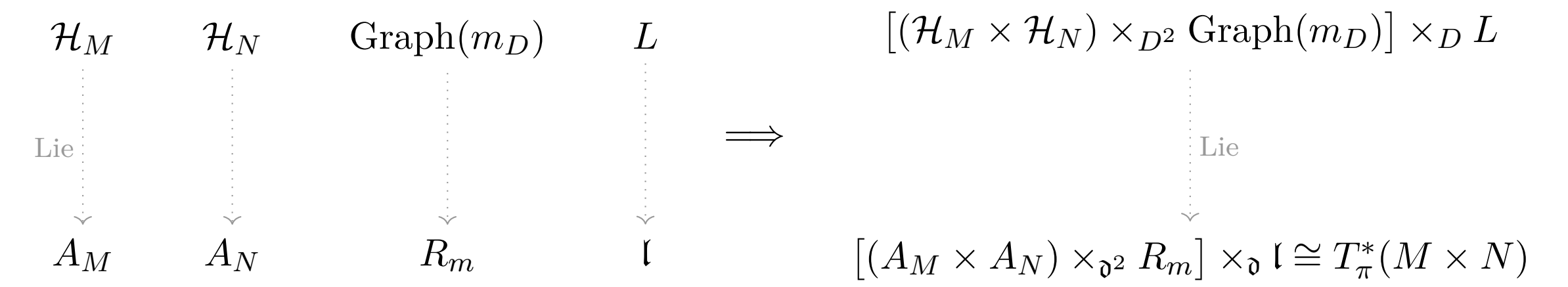
\begin{figure}[htbp]
\centering
\(
\begin{tikzcd}[column sep=1.4em, row sep=4em]
  \mathcal H_M
    \arrow[d, dotted, color=black!40, "\operatorname{Lie}"']
  & \mathcal H_N
    \arrow[d, dotted, color=black!40]
  & \operatorname{Graph}(m_D)
    \arrow[d, dotted, color=black!40]
  & L
    \arrow[d, dotted, color=black!40]
  \\
  A_M & A_N & R_m & \mathfrak l
\end{tikzcd}
\quad\Longrightarrow\quad
\begin{tikzcd}[row sep=4em]
  \bigl[
    (\mathcal H_M\times\mathcal H_N)
    \times_{D^2}\operatorname{Graph}(m_D)
  \bigr]\times_D L
    \arrow[d, dotted, color=black!40, "\operatorname{Lie}"]
  \\
  \bigl[
    (A_M\times A_N)\times_{\mathfrak d^2}R_m
  \bigr]\times_{\mathfrak d}\mathfrak l
  \cong T^*_{\pi}(M\times N)
\end{tikzcd}
\)\label{fig:intro}
\caption{Integration of Poisson structures coming from $r$-matrices, the arrows labelled by Lie indicate differentiation}
\end{figure}

\paragraph[Main results]{Main results.} In this note, we consider two further constructions of Poisson manifolds that may be put into the framework of such a theory of multiplicative $D$-valued moment maps. These are both associated with the Lie bialgebras $(\g,\g^*)$ determined by quasi-triangular $r$-matrices: these are elements in $\g \otimes \g$ whose respective symmetric and skew-symmetric parts $s$ and $\Lambda$ satisfy that $s$ is $\g$-invariant and both are linked by the classical Yang-Baxter equation. The key fact underlying this work is that $s$ determines a distinguished lagrangian Lie subalgebra 
\[ R_m \subset \db \times \db \times \overline{\db}, \]
where the overline indicates the opposite pairing on the corresponding factor. As we shall see, this lagrangian Lie subalgebra is responsible, firstly, for the
relationship between Poisson brackets and generalized dynamical
\(r\)-matrices pointed out in \cite{Xu}. Indeed, a generalized dynamical $r$-matrix over a Poisson manifold $(N,\pi_N)$ may also be encoded as a Lie algebroid 
\[A_N=T^*N\oplus (N\times \g^*)\] 
equipped with a distinguished Lie algebroid morphism to $\db$ \cite{DiracXu,Coupled}. Therefore, one can combine this Lie algebroid with the one coming from a Poisson action of $(\g,\g^*)$ on a Poisson manifold $(M,\pi_M)$ to produce a new Poisson bracket, see Proposition \ref{dynprop}. This construction gives rise to interesting examples, as pointed out already in \cite{MS,Xu}, it includes suitably gauge-fixed Poisson structures on moduli spaces of flat bundles over surfaces. Secondly, we consider the mixed
product Poisson structures of \cite{LuMouquin} and we show that they are also determined by suitable fibred products of Lie algebroids involving $R_m$. The main observation of this work is that then both constructions may be readily promoted to global ones in terms of $D$-valued moment maps. 
\begin{thm}\label{mainthm}
Let \(s+\Lambda\) be a quasitriangular \(r\)-matrix on \(\g\), and let
\(A_M\) be the Lie algebroid associated with a Poisson action of the
resulting Lie bialgebra on \((M,\pi_M)\). Suppose that \(A_M\) and the
Lie algebroid \(A_N\) in either of the following cases are integrable.
\begin{enumerate}
\item If \(A_N\) is associated with a generalized dynamical \(r\)-matrix
on \((N,\pi_N)\) with the twist determined by \(s\), then the coupled Poisson
structure \eqref{dyn} on \(M\times N\) is integrable.
\item If \(A_N\) is associated with a second Poisson action of the
same Lie bialgebra on \((N,\pi_N)\), then the mixed product
\eqref{mixed} on \(M\times N\) is integrable.
\end{enumerate}
\end{thm}
Moreover, in both cases there is an explicit construction of the integration given as follows. Let \(\Phi_M:\HH_M\to D\) and \(\Phi_N:\HH_N\to D\) be the Lie groupoid morphisms from the respective 
source-simply connected integrations of $A_M$ and $A_N$ determined by the above mentioned Lie algebroid morphisms. This construction relies crucially on a Lie \(2\)-group structure on $D$ integrating the double, which is guaranteed to exist if $D$ is simply-connected. In fact, if $D$ is simply-connected, it turns out that the Lie algebra $R_m$ may be promoted to a Lagrangian Lie subgroup of $D^3$ which is also the graph of a groupoid multiplication \(m_D\) on \(D\rightrightarrows G\), see \sect{global}.
In the dynamical case, let \(J_L:L\to D\) be the unit map \(G\to D\).
In the mixed case, let \(L=G^*\) be the simply connected dual group and
let \(J_L\) integrate the inclusion of its Lie algebra in \(\db\).
In either case, the arrows of a symplectic groupoid integrating the
indicated bivector are the triples
\((h_M,h_N,\ell)\in\HH_M\times\HH_N\times L\) for which
\(\Phi_M(h_M)\) and \(\Phi_N(h_N)\) are composable and
\[
 m_D\bigl(\Phi_M(h_M),\Phi_N(h_N)\bigr)=J_L(\ell).
\]
If \(\beta_M,\beta_N\) are the lagrangian two-forms on \(\HH_M,\HH_N\), the symplectic
form is the pullback of \(\pr_1^*\beta_M+\pr_2^*\beta_N\) under the
projection \((h_M,h_N,\ell)\mapsto(h_M,h_N)\). Schematically, our integration takes the form outlined in Figure~\ref{fig:intro}.

In order to illustrate this result with concrete examples, we revisit a simple construction of dynamical Poisson groupoids studied in \cite{LiParmentier}, and we recover the integration of configuration spaces of flags by means of decorated moduli spaces conducted in \cite{Morita}. For comparison, note that Mouquin constructed local symplectic groupoids over mixed product Poisson manifolds
\cite[Thm.~5.4 and Rem.~5.6]{MouquinLocal}. Our construction refines that observation by giving a corresponding global
symplectic groupoid when the two extended Lie algebroids encoding the Poisson action are integrable,
without assuming completeness of the dressing actions.

This paper is organized as follows. In \sect{inf}, we carry out the description of both procedures in terms of Lie algebroids and their fibred products.
In \sect{global}, we lift both constructions to the level of Lie groupoids, and we explain how they may give rise to new symplectic groupoids, by using the fact that a transverse intersection of
two 2-shifted lagrangians in a 2-shifted symplectic target is
1-shifted symplectic \cite[Thm.~2.9]{PTVV},
(see \cite[Thm.~3.23]{PartI} for a differential-geometric version of this result). We illustrate our results with explicit symplectic groupoids in both cases.
\paragraph*{Acknowledgments.} The author thanks H. Bursztyn, M. Cueca and M. Mol for stimulating conversations related to this work.
\section[Poisson brackets from generalized dynamical r-matrices and mixed products]{Poisson brackets from generalized dynamical
\(r\)-matrices and mixed products}\label{inf}
The central object in this work is given by a Lie algebra equipped with a quasi-triangular $r$-matrix. Such objects give rise to a special kind of Lie bialgebra which includes most important examples and, most crucially for us, allows us to perform a classical version of fusion for quantum group representations. This fusion underlies the mixed products of \cite{LuMouquin} and, as we shall see, the construction of Poisson structures from generalized dynamical \(r\)-matrices in \cite{Xu}.

Let \(\g\) be a Lie algebra equipped with a quasitriangular $r$-matrix 
\(s+\Lambda\in\g\otimes\g\). Here \(s\in(S^2\g)^{\g}\) is its
invariant symmetric part, which is possibly degenerate, \(\Lambda\in\wedge^2\g\) is its skew
part, and \(s+\Lambda\) solves the classical Yang--Baxter equation:
\[\tfrac12[\Lambda,\Lambda]=-\chi, \] 
where we write \(S=s^\sharp:\g^*\to\g\) for
contraction in the first factor and $\chi\in \wedge^3\g$ is given by
\(\chi(\xi,\eta,\zeta)=\langle\zeta,[S\xi,S\eta]\rangle\). 

We may use $\chi$ to define a Lie algebra structure on
\(\db=\g\oplus\g^*\), the double of \((\g,0,\chi)\). The Lie bracket is given by 
\[  [(u,\xi),(v,\eta)]= ([u,v]-i_\xi i_\eta \chi,\ad^*_u\eta -\ad^*_v\xi), \]
where $\ad_w^* \zeta(z)=-\zeta([w,z])$ and its pairing is
\[\langle(u,\xi),(v,\eta)\rangle=\xi(v)+\eta(u). \]
The dual Lie algebra of the Lie bialgebra defined by \(s+\Lambda\) sits in
this splitting as the lagrangian Lie subalgebra
\begin{equation}
       \mathfrak c=\{(\Lambda^\sharp\xi,\xi):\xi\in\g^*\}.
\label{dual lie algebra} \end{equation}
\paragraph[The fusion relation on the double]{The fusion relation on $\db$.} The double \(\db\) admits a groupoid structure
\(\db\rightrightarrows\g\) compatible with its Lie algebra structure.
We describe it following the conventions of \cite[Ex.~1.8, (1.9)]{LBS}. Its source, target and unit are
\[
 p_-(u,\xi)=u-S\xi,\qquad p_+(u,\xi)=u+S\xi,\qquad e(u)=(u,0).
\]
Two arrows \(d_1,d_2\) are composable when
\(p_-(d_1)=p_+(d_2)\). Their multiplication, in this order, is
\[m(d_1,d_2)=d_1+d_2-e(p_-(d_1)). \] 
If we let \(R_m\) be the graph of this multiplication, we have that it is a lagrangian Lie subalgebra of
\(\db\times\db\times\overline{\db}\), where the bar indicates that the
pairing is taken with an opposite sign. Write \(p:R_m\to\db^2\) for projection to the pair of composable elements and
\(q:R_m\to\db\) for the projection to the output of the multiplication.

For example, if \(s\) is nondegenerate, the map
\((p_+,p_-):(u,\xi)\mapsto(u+S\xi,u-S\xi)\) identifies
\(\db\rightrightarrows\g\) with the pair groupoid
\(\g\times\g\rightrightarrows\g\). In terms of its natural coordinates, we have that \((x,y)\)
is an arrow from \(y\) to \(x\), the unit at \(x\) is \((x,x)\), and the multiplication is
\((x,y)\cdot(y,z)=(x,z)\).

The Lie algebra $R_m$ gives us the ability to carry out the following operation. Suppose \(A_i\to X_i\), \(i=1,2\), are Lie algebroids equipped with Lie algebroid morphisms
\(\varphi_i:A_i\to\db\). For a lagrangian Lie 
subalgebra \(\mathfrak l\) with inclusion \(j:\mathfrak l\to\db\),
consider the following two fibred products of Lie algebroids over \(X_1\times X_2\) which we assume to be transverse:
\begin{equation}\label{squares}
\begin{tikzcd}[column sep=large,row sep=small]
 A_{\fus}\arrow[r]\arrow[d] & R_m\arrow[d,"p"]\\
 A_1\times A_2\arrow[r,"\varphi_1\times\varphi_2"'] & \db^2
\end{tikzcd}
\hspace{2em}
\begin{tikzcd}[column sep=large,row sep=small]
 A_{\mathfrak l}\arrow[r]\arrow[d] &\mathfrak l\arrow[d,"j"]\\
 A_{\fus}\arrow[r,"\varphi_{\fus}"'] &\db .
\end{tikzcd}
\end{equation}
Here \(\varphi_{\fus}\) is the upper horizontal map in the first square
followed by \(q\). The first square serves to multiply the target elements of $\varphi_1\times\varphi_2$; on the other hand,
the second square further constrains the multiplication of a pair of elements to lie in \(\mathfrak l\).
In the two examples below, we shall study special situations of such fibred products in which we may identify \(A_{\mathfrak l}\) with the graph of a Poisson tensor.  

\paragraph[Generalized dynamical \texorpdfstring{\(r\)}{r}-matrices and Poisson brackets]{Generalized dynamical \(r\)-matrices and Poisson brackets.}
In mathematical physics, classical dynamical \(r\)-matrices occur as semiclassical limits of
dynamical quantum \(R\)-matrices, including Felder's elliptic examples
\cite{EV}. This led to an explicit relationship with Lie theory and Poisson geometry.
In \cite{Coupled,LiuXu} the dynamical equation was generalized to the case of a Poisson manifold as base, which is the situation relevant for us.

Let \((N,\pi_N)\) be Poisson. Choose a basis \((e_i)_i \) of \(\g\).
A generalized dynamical datum with the twist fixed above by $\chi$ consists of
a tensor \(\vartheta=\sum_i X_i\otimes e_i\), with \(X_i\) vector
fields on \(N\), and a smooth map \(r:N\to\wedge^2\g\).
Write \(b=\vartheta^\sharp:T^*N\to\g\), then we have that
\(b\alpha=\sum_i\alpha(X_i)e_i\). Consider the Lie algebroid
\(E=TN\oplus(N\times\g)\) with anchor the projection
to \(TN\), and its bracket on sections defined by
\[
 [(X,u),(Y,v)]_E
 =\bigl([X,Y],\,dv(X)-du(Y)+[u,v]_{\g}\bigr),
\]
where \(u,v:N\to\g\) and \([u,v]_{\g}\) is the pointwise Lie
bracket \cite[(28)]{Xu}, with vector fields acting on
\(\g\)-valued functions by differentiation. We use the same notation for the induced
Schouten bracket on \(\Gamma(\wedge^\bullet E)\). Put
\(\mathcal R=\pi_N+\sum_i X_i\wedge e_i+r\). The compatibility
condition for these data is
\[\tfrac12[\mathcal R,\mathcal R]_E=-\chi. \]
In particular, it follows that \(b:T^*_{\pi_N}N\to\g\) is a Lie algebroid morphism
\cite[Def.~4.4 and Thm.~4.5]{LiuXu}.
Moreover, the dynamical datum defines a Lie algebroid
\(A_N=T^*N\oplus(N\times\g^*)\) and a Lie algebroid morphism
\[
 \varphi_N:A_N\longrightarrow\db,\qquad
 \varphi_N(\alpha,\xi)=(b\alpha+r^\sharp\xi,\xi).
\]

Similarly, we may construct a Lie algebroid from a Poisson action.
Let \((M,\pi_M)\) be a Poisson manifold with a Poisson action of the
Lie bialgebra with cobracket \(\delta(u)=[u,\Lambda]\) determined by \(s+\Lambda\). We denote by
\(\rho:\g\to\mathfrak X(M)\) the corresponding
action. In the convention of
\cite[Ex.~3.19]{PartI}, the Poisson action condition is
\[
 \mathcal L_{\rho(u)}\pi_M=-\rho(\delta(u)),\qquad u\in\g.
\]
The associated Lie algebroid \(A_M=T^*M\oplus(M\times\g)\) has anchor
\(\mathsf a_M(\beta,u)=\pi_M^\sharp\beta+\rho(u)\).
In the splitting \(\db=\g\oplus\g^*\) used above, we obtain a Lie algebroid
morphism to $\db$ given by
\[
 \varphi_M:A_M\longrightarrow\db,\qquad
 \varphi_M(\beta,u)=(u-\Lambda^\sharp\rho^*\beta,-\rho^*\beta),
\]
where we denote by
\(\rho^*:T^*M\to\g^*\) the fibrewise transpose of the action map. We may describe both Lie algebroids uniformly by embedding them as Dirac
structures. Denote by \(\TT X=TX\oplus T^*X\) the standard Courant
algebroid of a manifold \(X\). The embeddings of \(A_M\) into
\(\TT M\times\db\) and \(A_N\) into \(\TT N\times\db\) take the following
form, see \cite[(43), (45)]{PartI}:
\begin{align}
 (\beta,u)&\longmapsto
 \bigl((\pi_M^\sharp\beta+\rho(u),\beta),\varphi_M(\beta,u)\bigr),
 \notag\\
 (\alpha,\xi)&\longmapsto
 \bigl((\pi_N^\sharp\alpha-b^*\xi,\alpha),\varphi_N(\alpha,\xi)\bigr).
 \label{embeddings}
\end{align}
The Lie brackets on $A_M$ and $A_N$ are then obtained by restriction of the background Courant brackets, and the anchors are their tangent projections. For generalized
dynamical \(r\)-matrices, this Dirac description goes back to
\cite{DiracXu}, see also \cite[Ex.~3.20]{PartI}. Applying construction \eqref{squares} to these Lie algebroids, we recover the Poisson structure of \cite[Thm. 2.9]{Xu}. For a general interpretation of such fibred products of Lie algebroids in terms of Dirac reduction see \cite{PartI}.

\begin{prop}\label{dynprop}
Take \((A_1,A_2)=(A_M,A_N)\) as in \eqref{embeddings} and
\(\mathfrak l=e(\g)\) in \eqref{squares}. The resulting Lie algebroid is the cotangent
algebroid of the Poisson bivector
\begin{equation}\label{dyn}
 \Pi=Q+\pi_N-\sum_i\rho(e_i)\wedge X_i+\rho(r)
 \qquad\text{on } M\times N .
\end{equation}
Here \(Q=\pi_M+\rho(\Lambda)\).
\end{prop}
\begin{proof} 
Put \(d_M=\varphi_M(\beta,u)\) and \(d_N=\varphi_N(\alpha,\xi)\).
The first fibred product in \eqref{squares} imposes \(p_-(d_M)=p_+(d_N)\).
The second fibred product requires \(m(d_M,d_N)\) to be a unit, so the two arrows $d_M$ and $d_N$ are
inverses, which implies that their \(\g^*\)-components sum to zero. Therefore, we have that
\[
 \xi=\rho^*\beta,\qquad u=b\alpha+(r+\Lambda)^\sharp\rho^*\beta .
\]
Consequently, each covector \((\beta,\alpha)\in T^*(M\times N)\) has a unique lift
through the two fibred products in\eqref{squares}. Its tangent projection is then given by the anchor
\[
 \bigl(\pi_M^\sharp\beta+\rho\bigl(b\alpha+(r+\Lambda)^\sharp\rho^*\beta\bigr),\;
 \pi_N^\sharp\alpha-b^*\rho^*\beta\bigr)
       =\Pi^\sharp(\beta,\alpha).
\]
The resulting Lie algebroid is therefore isomorphic to \(\Gr(\Pi^\sharp)\). Indeed, projection to \((\beta,\alpha)\) gives the Lie algebroid isomorphism
\[
 A_{e(\g)}
 =\bigl((A_M\times A_N)\times_{\db^2}R_m\bigr)\times_{\db}e(\g)
 \ \cong\ T^*_\Pi(M\times N).
\]
The flip \(N\times M\to M\times N\)
takes the Poisson structure in \cite[Thm.~2.9, (14)]{Xu} to \(\Pi\). \end{proof}

\paragraph[Example. A gauge-fixed Chern--Simons bracket]{Example. A gauge-fixed Chern--Simons bracket.}  We shall illustrate how Proposition \ref{dynprop} gives rise to interesting examples associated to the Fock-Rosly Poisson structure on moduli spaces of flat bundles over surfaces, following \cite{MS,Xu}. Let \(G=\mathrm{SO}_+(2,1)\ltimes\mathbb R^3\), and consider an oriented connected surface $\Sigma$
of genus \(g\) with \(n\geq3\) punctures. Choose a conjugacy class
\(C_i\subset G\) for the holonomy at puncture \(i\). On the regular
locus of \(P=C_1\times C_2\), take a local slice \(N\subset P\) where simultaneous
conjugation is locally free. The remaining holonomies of a representation of $\pi_1(\Sigma)$ belong then to
\[M=C_3\times\cdots\times C_n\times G^{2g}.  \]
Let \(\rho_P\) denote infinitesimal conjugation on \(P\), and let
\(p_N:TP|_N\to TN\) be the projection along the conjugation orbits. We have that $P$ inherits a canonical bivector field \(Q_P\) given as follows. Choose a
basis \((e_a)\) of \(\mathfrak g\) and write the symmetric part of the
quasitriangular tensor as \(s=\sum_{a,b}s^{ab}e_a\otimes e_b\), where \(2s\) is the inverse tensor of the Chern–Simons pairing.
For \(u\in\mathfrak g\), let \(u^L\) and \(u^R\) be its respective left and right invariant extensions. In what follows, a superscript \((i)\) denotes the lift to the \(i\)-th factor of
\(G\times G\). Thus \(\rho_i(u)=(u^L-u^R)^{(i)}\) generates the action by conjugation on that factor. The canonical (quasi-Poisson) bivector field on $P$ is
\[
  Q_P=
  \left.
  \sum_{a,b}s^{ab}\Bigl(
    (e_a^R)^{(1)}\wedge(e_b^L)^{(1)}
    +(e_a^R)^{(2)}\wedge(e_b^L)^{(2)}
    +\rho_1(e_a)\wedge\rho_2(e_b)
  \Bigr)\right|_{C_1\times C_2};
\]
see \cite[Example~2.11, Eq.~(15)]{Xu}. Now we obtain a Poisson structure on the slice via projection
\(\pi_N=(\wedge^2 p_N)(Q_P|_N)\), for the Dirac gauge fixing of \cite{MS} we have \(\pi_N=0\) \cite[Thm. 3.6]{Xu}. Furthermore, the mixed and orbit components of $Q_P$
determine \(\vartheta=\sum_i X_i\otimes e_i\) and \(r\) through the decomposition described in \cite[Thms.~2.7 and 3.4]{Xu}:
\[
 Q_P|_N=\pi_N-\sum_i X_i\wedge\rho_P(e_i)+\rho_P(r).
\]
We have then that $\pi_N$, $\vartheta$ and $r$ give us a generalized dynamical $r$-matrix on $N$.

The Poisson structure $\pi_M$ on the space of remaining holonomies
\(M=C_3\times\cdots\times C_n\times G^{2g}\) is given by the Fock-Rosly construction, which is determined by an $r$-matrix of the form \(s-\Lambda\), see \cite{Mouquin}. Hence, writing
\(\vartheta=\sum_i X_i\otimes e_i\), equation~\eqref{dyn} gives
the gauge-fixed Poisson bivector on \(M\times N\),
\[
  \Pi=\pi_M+\rho(\Lambda)+\rho(r)-\sum_i\rho(e_i)\wedge X_i ;
\]
here \(r\) is evaluated at the point of \(N\) and $\rho$ is the infinitesimal conjugation action on $M$
\cite[Thms.~3.4 and 3.6]{Xu}. We leave the global study of this example to \cite{PartII}.

\paragraph[The mixed product of two Poisson spaces]{The mixed product of two Poisson spaces.}
Write \(r=s+\Lambda=\sum_a x_a\otimes y_a\) for the fixed
quasitriangular tensor. Let \((M_1,\pi_1)=(M,\pi_M)\) and
\((M_2,\pi_2)=(N,\pi_N)\) both be left Poisson spaces for \(r\),
with infinitesimal actions \(\rho_1,\rho_2\). In the left-action convention of \cite[\S6.3]{LuMouquin}, the mixed
product of these Poisson structures is
\begin{equation}\label{mixed}
 \pi_{\mathrm{mix}}=\pi_1+\pi_2+
       \sum_a\rho_1(y_a)\wedge\rho_2(x_a).
\end{equation}

Consider the corresponding Lie algebroid for each factor \(A_i=T^*M_i\oplus(M_i\times\g)\).
For each \(a_i=(\alpha_i,v_i)\in A_i\), we write
\(\xi_i=-\rho_i^*\alpha_i\), where \(\rho_i^*\) is the transpose of
the action. In the splitting \(\db=\g\oplus\g^*\) used at the beginning of this section, the Dirac embedding of \cite[(43)]{PartI} is given by
\[
 a_i\longmapsto
 \big((\pi_i^\sharp\alpha_i+\rho_i(v_i),\alpha_i),
          (v_i+\Lambda^\sharp\xi_i,\xi_i)\big)
 \quad\text{in }\TT M_i\times\db .
\]
We abbreviate the component in \(\db\) by \(d_i=\varphi_i(a_i)\) in what follows.
\begin{prop}\label{mixprop}
Applying the two fibred products in \eqref{squares} successively to
\(A_1,A_2\) with \(\mathfrak l=\mathfrak c\) as in \eqref{dual lie algebra} determines the cotangent Lie
algebroid of the mixed product \eqref{mixed}:
\begin{equation}\label{mixsquares}
\begin{tikzcd}[column sep=large,row sep=large]
 F\arrow[r,"\pr_2"]\arrow[d,"\pr_1"'] &
 A_2\arrow[d,"p_+\circ\varphi_2"]\\
 A_1\arrow[r,"p_-\circ\varphi_1"'] & \g
\end{tikzcd}
\hspace{2em}
\begin{tikzcd}[column sep=large,row sep=large]
 C\arrow[r]\arrow[d,hook] &
 \mathfrak c\arrow[d,hook,"j"]\\
 F\arrow[r,"\varphi_F"'] & \db ;
\end{tikzcd}
\end{equation}
where the projection to the cotangent bundle factor identifies \(C\) with
\(T^*_{\pi_{\mathrm{mix}}}(M_1\times M_2)\).
\end{prop}
\begin{proof} Explicitly, we have that \(F\)
denotes the Lie algebroid formed by pairs \((a_1,a_2)\) satisfying
\(p_-(d_1)=p_+(d_2)\), this means that: 
\begin{equation*}
 v_2+(\Lambda^\sharp+S)\xi_2=v_1+(\Lambda^\sharp-S)\xi_1 .
\label{composability}\end{equation*}
The map $\varphi_F:F\to \db$ is then given by \(\varphi_F(a_1,a_2)=m(d_1,d_2)\). If we let \(j:\mathfrak c\hookrightarrow\db\)
be the inclusion, we have that \(C=F\times_{\db}\mathfrak c\). The condition defining the right square above may be written using the formula
\[
m(d_1,d_2)=\big(v+\Lambda^\sharp(\xi_1+\xi_2),\,\xi_1+\xi_2\big).
\]
where $v=v_2-\Lambda^\sharp\xi_1+S\xi_1$. Since \(\mathfrak c\) is the graph of \(\Lambda^\sharp\), the right
square sets \(v=0\). Thus we obtain that
\[
 v_1=(\Lambda^\sharp+S)\xi_2,\qquad
 v_2=(\Lambda^\sharp-S)\xi_1 .
\]
Since each covector \((\alpha_1,\alpha_2)\)
determines a unique element of \(C\), we may substitute the two values of
\(v_i\) in the tangent entries of the Dirac embeddings to obtain that
\[
 \left(
 \pi_1^\sharp\alpha_1-\sum_a\alpha_2(\rho_2(x_a))\rho_1(y_a),\;
 \pi_2^\sharp\alpha_2+\sum_a\alpha_1(\rho_1(y_a))\rho_2(x_a)
 \right).
\]
These are exactly the two components of
\(\pi_{\mathrm{mix}}^\sharp(\alpha_1,\alpha_2)\).
Therefore, we may identify \(C\) with
\(\Gr(\pi_{\mathrm{mix}}^\sharp)\).
\end{proof}
\begin{rema} The key feature of $\pi_{\mathrm{mix}}$ is that the diagonal $\g$-action on $M_1\times M_2$ becomes a Poisson action with respect to it \cite{LuMouquin} and the full Lie algebroid which encodes this action is precisely the fibred product $F$ above.    
\end{rema}
\section{The global picture: symplectic groupoids from fibred products}\label{global}

To integrate the two pullback diagrams in \eqref{squares}, we first integrate
\(\db\rightrightarrows\g\). To do so we use the standard integration of crossed modules of Lie algebras. Introduce the new coordinates
\((x,\xi)=(u-S\xi,\xi)\) in which the source map is \(x\) and the target is
\(x+2S\xi\). Expressing the coadjoint action as
\((\ad^*_x\xi)(y)=-\xi([x,y])\), we obtain the
crossed module of Lie algebras 
\[\partial=2S:\mathfrak k\to\g,\] 
where \(\mathfrak k\) is the vector space \(\g^*\) with bracket
\([\xi,\eta]_{\mathfrak k}=\ad^*_{\partial\xi}\eta\), and \(\g\)
acts on it by the coadjoint action above. In this expression, the pairing becomes
\[
 \xi(y)+\eta(x)+2s(\xi,\eta).
\]

Let \(K,G\) be the simply connected groups integrating
\(\mathfrak k,\g\), and use \(\partial:K\to G\) also for the
integrated homomorphism. The semidirect product Lie group \(D=K\rtimes G\) is in addition a
groupoid over \(G\), with source \(g\), target \(\partial(k)g\)
for \((k,g)\), and unit \(\epsilon(g)=(1,g)\). The presence of this groupoid structure makes $D$ into a {\em Lie 2-group}, a groupoid object in the category of Lie groups. In particular, if \(S\) is invertible,
the target and source identify this groupoid with the pair groupoid
\(G\times G\rightrightarrows G\). 

Equip \(D\) with its standard 2-shifted symplectic form, see
\cite[Ex.~3.1]{PartI} for instance. We let \(\theta^L,\theta^R\) be the respective left and
right Maurer--Cartan forms on \(D\) and define
\[
 \Omega_D=-\tfrac12\langle\pr_1^*\theta^L,\pr_2^*\theta^R\rangle\in \Omega^2(D\times D),
 \qquad
 \Theta_D=-\tfrac1{12}\langle\theta^L,[\theta^L,\theta^L]\rangle \in \Omega^3(D).
\]
We say that $\Omega_D$ is the Polyakov-Wiegmann 2-form and $\Theta_D$ is the Cartan 3-form. Let \(R_{m_D}\subset D^2\times\overline D\) be the graph of the groupoid multiplication \(m_D\), written as a relation from $D^2$ to $D$; here the bar reverses
the sign of the shifted symplectic form above. We have that $R_{m_D}$ is a Lie group with Lie algebra \(R_m\). In the terminology of \cite{PartI}, this subgroup is 2-shifted lagrangian with zero two-form: we may readily verify its Lie
algebra is lagrangian (see e.g. \cite[(1.9)]{LBS}), and the differential forms induced by $ \Omega_D $ and $\Theta_D$ on the nerve of $D^2\times \overline{D}$ restrict to zero on (the nerve of) $R_{m_D}$ by isotropy. Here the overline indicates a sign reversal in such differential forms for the corresponding factor. 

The second idea needed to carry out the integration of \eqref{squares} is that of lagrangian structures and their corresponding notion of composition, which we summarize from \cite{PartI}. Let \(\Phi:\HH\to D\) be a Lie groupoid morphism, where
\(\HH\rightrightarrows X\). Write \(A_{\HH}\) for the Lie algebroid of \(\HH\),
\(\mathsf a:A_{\HH}\to TX\) for its anchor, and
\(\varphi=\Lie(\Phi):A_{\HH}\to\db\) for the induced morphism.
A two-form \(\beta\in\Omega^2(\HH)\) determines a bundle map
\(\lambda_\beta:A_{\HH}\to T^*X\) by
\[
  \lambda_\beta(a)(v)=\beta_{1_x}\bigl(a,T_x\varepsilon(v)\bigr),
  \qquad a\in(A_{\HH})_x,\quad v\in T_xX,
\]
where \(\varepsilon:X\to\HH\) is the unit map. We say that \(\beta\) defines a
{\em \(2\)-shifted lagrangian structure} on \(\Phi\) if
\[
  \pr_1^*\beta+\pr_2^*\beta-m_{\HH}^*\beta
    =(\Phi\times\Phi)^*\Omega_D,
  \qquad
  -d\beta=\Phi^*\Theta_D,
\]
and the bundle map
\[
  (\mathsf a,\lambda_\beta,\varphi):
  A_{\HH}\longrightarrow \TT X\oplus(X\times\db)
\]
is fibrewise injective with lagrangian image
\cite[Prop.~3.9]{PartI}. The first identity above takes place on the space of
composable pairs in \(\HH\), with projections \(\pr_1,\pr_2\)
and multiplication \(m_{\HH}\). When \(D=\{1\}\), these conditions
say that \((\HH,\beta)\) is a quasi-symplectic groupoid with zero
background three-form. One of the main results in \cite{PartI} is that integrable Dirac structures in $\TT X\oplus(X\times\db)$ give rise upon integration to such lagrangian structures.

We recall how to compose two such lagrangian structures, in analogy with usual lagrangian correspondences in classical symplectic geometry. Let \(D_1,D_2,D_3\)
be Lie groups equipped with \(2\)-shifted symplectic forms as before, and write \(\overline D_i\)
for \(D_i\) equipped with the opposite shifted symplectic form. Suppose that
\[
  \Phi_{12}=(f_1,f_2):\HH_{12}\longrightarrow\overline D_1\times D_2,
  \qquad
  \Phi_{23}=(g_2,g_3):\HH_{23}\longrightarrow\overline D_2\times D_3
\]
carry lagrangian structures \(\beta_{12}\) and \(\beta_{23}\),
respectively. If \(f_2\) and \(g_2\) are transverse, let
\(\HH_{13}=\HH_{12}\times_{D_2}\HH_{23}\). Then
\[
  \Phi_{13}=(f_1,g_3):
  \HH_{13}\longrightarrow\overline D_1\times D_3
\]
carries the lagrangian structure
\[
  \beta_{13}=\text{pr}_1^*\beta_{12}+\text{pr}_2^*\beta_{23},
\]
see \cite[Thm.~3.23]{PartI}. 
\begin{exa}[\cite{PartI}]\label{ex:poiact} A recurrent example of this construction is the following. Suppose that the Lie algebra inclusion $\g^*\subset \db$ integrates to a Lie group morphism $G^*\to D$. Then if $\beta$ is a lagrangian structure on $\Phi: \HH\to D$ integrating the Lie algebroid of a Poisson action of $(\g,\g^*)$ on $(M,\pi_M)$, we have that $\HH\times_D G^*$ equipped with the restriction of $\pr_1^*\beta$ is a symplectic groupoid integrating $\pi_M$.     
\end{exa}

We may apply this composition rule twice to integrate the two
fibred products in \eqref{squares}. Let
\(\Phi_i:\HH_i\to D\) be the integrating morphisms with two-forms
\(\beta_i\), and let \(J_L:L\to D\) integrate the inclusion
\(\mathfrak l\hookrightarrow\db\). Write \(p:R_{m_D}\to D\times D\)
and \(q:R_{m_D}\to D\) for the input and output projections of the
multiplication graph. The three correspondences are
\[
\begin{tikzcd}[column sep=small,row sep=large]
 & \HH_1\times\HH_2
   \arrow[dl]\arrow[dr,"{\Phi_1\times\Phi_2}"]
 && R_{m_D}\arrow[dl,"p"']\arrow[dr,"q"]
 && L\arrow[dl,"J_L"']\arrow[dr] & \\
 \{1\} && D\times D && D && \{1\}.
\end{tikzcd}
\]
Since \(R_{m_D}\) and \(L\) carry lagrangian structures given by the zero 2-form, it follows that the
resulting groupoid
\[
  \mathcal K_L=
  \bigl((\HH_1\times\HH_2)\times_{D\times D}R_{m_D}\bigr)\times_D L
\]
is equipped with the symplectic form \(\omega=\pr_1^*\beta_1+\pr_2^*\beta_2\), where
\(\pr_i:\mathcal K_L\to\HH_i\) are the projections.
\begin{thm}\label{globalprop}
Suppose that the Lie algebroids $A_1$ and $A_2$ in Proposition~\ref{dynprop}
or Proposition~\ref{mixprop} are integrable. Then the corresponding
Poisson manifold is integrable. In both cases a symplectic groupoid
is obtained by replacing each Lie algebroid in
\eqref{squares} by its groupoid counterpart.
\end{thm}
\begin{proof}
For each Lie algebroid \(A_i\), take its source-simply connected integration
\(\HH_i\rightrightarrows X_i\). The corresponding Dirac embedding determines a lagrangian structure \(\beta_i\) on the Lie groupoid morphism
\(\Phi_i:\HH_i\to D\) integrating \(\varphi_i\)
\cite[Cor.~3.13]{PartI}. In the case of the Lie algebroids coming from Poisson actions, \(\Phi_i\) is, in addition, a
multiplicative $D$-valued moment map as in \cite[Thm.~5.3]{PartI}. Let \(P:R_{m_D}\to D^2\) be the projection to the composable pair and let
\(\Phi_{\fus}\) be the projection to the multiplication after the first pullback below.
Denoting by \(J_L:L\to D\) the integration of the inclusion
\(\mathfrak l\subset \db\) used in \eqref{squares}, the relevant fibred products are
\begin{equation}\label{globalsquares}
\begin{tikzcd}[column sep=large,row sep=small]
 \HH_{\fus}\arrow[r]\arrow[d] & R_{m_D}\arrow[d,"P"]\\
 \HH_1\times\HH_2\arrow[r,"\Phi_1\times\Phi_2"'] & D^2
\end{tikzcd}
\hspace{2em}
\begin{tikzcd}[column sep=large,row sep=small]
 \mathcal K_L\arrow[r]\arrow[d] & L\arrow[d,"J_L"]\\
 \HH_{\fus}\arrow[r,"\Phi_{\fus}"'] & D .
\end{tikzcd}
\end{equation}
For the Poisson structure coming from the dynamical $r$-matrix take \(J_L=\epsilon:G\to D\).
For the mixed product Poisson structure take \(L=G^*\), the simply connected group
integrating \(\mathfrak c\), with its integrating immersion \(J_L\).
Thus we obtain \(\mathcal K_L=\HH_{\fus}\times_D L\).

In Propositions~\ref{dynprop} and~\ref{mixprop}, each of the two
fibred product conditions has the expected codimension
\(\dim\g\).
This gives transversality on the source fibres at the units, and
right translation implies that transversality holds at every other groupoid element. Thus both
pullbacks in \eqref{globalsquares} are smooth Lie groupoids. By construction, differentiating \eqref{globalsquares} recovers \eqref{squares}:
\[ \Lie(R_{m_D})=R_m, \qquad \Lie(m_D)=m, \qquad \Lie(J_L)=j. \]
By \cite[Thms.~3.22--3.23, Cor.~3.24]{PartI}, the corresponding lagrangian
structures induce on \(\mathcal K_L\) the restriction of
\(\pr_1^*\beta_1+\pr_2^*\beta_2\) as discussed above and it integrates the graph of the Poisson structures identified in
Proposition~\ref{dynprop} or Proposition~\ref{mixprop}; hence this 2-form
is symplectic.
\end{proof}

\begin{cor}
Let \(\mathfrak h\subset\g\) be a Lie subalgebra and \(N\) a connected
open subset of \(\mathfrak h^*\) with its linear Poisson structure and
\(\pi_2(N)=0\). A classical dynamical \(r\)-matrix on \(N\) determines an
integrable Poisson structure \eqref{dyn} for every integrable input \(A_M\).
\end{cor}
\begin{proof}
Under the identification \(T^*N=N\times\mathfrak h\), we obtain that the map \(b\) is the inclusion
of \(\mathfrak h\) in \(\g\), hence its transpose is surjective, and so the
anchor of \(A_N\) in \eqref{embeddings} is surjective. The resulting
transitive Lie algebroid has trivial monodromy because \(\pi_2(N)=0\)
\cite{CF}. It follows that Theorem~\ref{globalprop} may be applied.
\end{proof}
For iterated mixed products, we need to integrate the full Lie algebroid
associated with the Poisson action at each stage.

\begin{cor}\label{cor:mix}
Suppose that the full Lie algebroids associated with the two Poisson
actions in Proposition~\ref{mixprop} are integrable. Then the first
fibred product in \eqref{globalsquares},
\[
  \HH_{\fus}
  = (\HH_1\times\HH_2)\times_{D\times D} R_{m_D},
\]
integrates the full Lie algebroid associated with the diagonal Poisson
action on their mixed product.
\end{cor}

\begin{proof}
Keeping \(v\in\g\) free in the proof of Proposition~\ref{mixprop}
adds \((\rho_1(v),\rho_2(v))\) to the cotangent anchor of the mixed
product and also gives the Dirac embedding of the diagonal action.
The first fibred product in \eqref{globalsquares} integrates this embedding.
\end{proof}

\paragraph[Example. The rational dynamical \texorpdfstring{\(r\)}{r}-matrix for \texorpdfstring{\(\mathfrak{sl}_2\)}{sl(2)}]{Example. The rational dynamical \(r\)-matrix for \(\mathfrak{sl}_2\).}
Choose a basis \(H,E_+,E_-\) of \(\g=\mathfrak{sl}_2(\mathbb R)\)
satisfying \([H,E_\pm]=\pm2E_\pm\) and \([E_+,E_-]=H\).
On the dynamical base \(N=U=\{\lambda\in(\mathbb RH)^*:q=\lambda(H)>0\}\)
take \(\pi_U=0\), \(s=\Lambda=0\) and
\[
 \vartheta=\partial_q\otimes H,\qquad
 r(q)=-q^{-1}E_+\wedge E_- .
\]
Indeed, for \(\mathcal R=\partial_q\wedge H+f(q)E_+\wedge E_-\),
the coefficient of \(H\wedge E_+\wedge E_-\) in
\(\tfrac12[\mathcal R,\mathcal R]_E\) is \(f^2-f'\), and all
other coefficients vanish. The choice \(f(q)=-1/q\) therefore
satisfies the dynamical equation in \sect{inf}.

For \(a\in\mathbb R\) and \(\xi\in\g^*\), write
\(\mathcal R^\sharp(a\,dq,\xi)=(c\,\partial_q,u)\). Then
\[
 c=-\xi(H),\qquad
 u=aH+q^{-1}\bigl(\xi(E_-)\,E_+-\xi(E_+)\,E_-\bigr).
\]
These equations determine \(a,\xi\) uniquely from \(c,u\).
Since \([\mathcal R,\mathcal R]_E=0\), the map
\(\mathcal R^\sharp\) identifies the Dirac algebroid \(A_U\) of
\eqref{embeddings} with the gauge algebroid \(TU\oplus(U\times\g)\)
\cite[\S2]{PingTriangular}. For a classical dynamical base, \(A_U\)
is the dual algebroid of the dynamical Poisson groupoid
\cite[Ex.~4.9]{PartI}. Li and Parmentier describe the dual groupoid
over a contractible base \cite[Thm.~3.2.4]{LiParmentier}.

Let \(G\) be the simply connected group with Lie algebra \(\g\).
The source-simply connected integration in this example is
\[\HH_U=U\times G\times U\rightrightarrows U. \] 
An arrow
\((y,g,x)\) has source \(x\) and target \(y\), and the multiplication is given by \((z,h,y)(y,g,x)=(z,hg,x)\). Each source fibre is thus of the form \(U\times G\).
Write \(\eta\in\g^*\) for the covector
with \(\eta(H)=1\) and \(\eta(E_\pm)=0\), and put \(\ell(q)=q\eta\).
Since \(s=0\), we can write the double with the group coordinate first
as \(D=G\ltimes\g^*\), where
\(\operatorname{Ad}_g^*\xi=\xi\circ\operatorname{Ad}_{g^{-1}}\).
Let \(\theta_G^R=dg\,g^{-1}\) and \(\theta_G^L=g^{-1}dg\) denote the
matrix Maurer--Cartan forms pulled back to \(G\) from
\(\mathrm{SL}_2(\mathbb R)\). The moment map and its lagrangian structure are
\begin{equation}\label{rational-integration}
 \begin{aligned}
 \Phi_U(y,g,x)&=\bigl(g,\ell(y)-\operatorname{Ad}_g^*\ell(x)\bigr),\\
 \beta_U&=\tfrac12\,dy\wedge\langle\eta,\theta_G^R\rangle
           +\tfrac12\,dx\wedge\langle\eta,\theta_G^L\rangle .
 \end{aligned}
\end{equation}
In order to verify that \(\Phi_U\) is a morphism, we may write it as the product of
\(A=(1,\ell(y))\), \(B=(g,0)\) and \(C=(1,-\ell(x))\), it thus follows that the two factors at the common endpoint cancel under composition. Since the subgroups \(G\) and \(\g^*\) are isotropic, the
Polyakov--Wiegmann identity for this factorization implies that
\(\beta_U=-(A,B)^*\Omega_D-(AB,C)^*\Omega_D\). This is the second
formula in \eqref{rational-integration} and it implies the two required 
identities for $\beta_U$.

At the unit \((q,1,q)\), the algebroid element \((c\,\partial_q,u)\)
is represented by \(v=(c,-u,0)\) in the right-invariant convention
of \cite[\S2]{PartI}. Using the same identification in the double,
we obtain
\[
 \Lie(\Phi_U)(c\,\partial_q,u)
   =\bigl(u,-c\eta-\ad_u^*\ell(q)\bigr)=(u,\xi),\qquad
 \beta_U\bigl(v,(Y,0,Y)\bigr)=\eta(u)Y=aY
 \quad(Y\in\mathbb R).
\]
Together with the anchor \(c\,\partial_q\), these are exactly the components of the embedding in
\eqref{embeddings} corresponding to the Dirac structure associated to a dynamical $r$-matrix. Therefore, we have that \((\HH_U,\Phi_U,\beta_U)\) integrates the
Dirac structure of \(r(q)\).

In order to apply Theorem \ref{globalprop}, we may take any integrable Lie algebroid \(A_M=T^*M\oplus (M\times \g)\) corresponding to a $G$-action on $M$ by Poisson automorphisms. Let \(\beta_M\) be the corresponding lagrangian structure and denote by \(\Phi_M(k)=(g_M(k),\nu_M(k))\) the components of its moment map, which takes values in \(D=G\ltimes\g^*\). The two conditions expressed by \eqref{globalsquares} now read
\[
   g_M(k)=g,\qquad \nu_M(k)+\ell(y)-\operatorname{Ad}_g^*\ell(x)=0 .
\]
It follows that the pairs \((k,(y,g,x))\) satisfying these equations constitute a groupoid
over \(M\times U\) integrating \eqref{dyn}, whose symplectic form is the
restriction of \(\pr_1^*\beta_M+\pr_2^*\beta_U\).

\paragraph[Example. Integration of the flag variety]{Integration of the flag variety.}
Let \(G\) be a connected complex semisimple group with opposite
Borel subgroups \(B=B_+\) and \(B_-\) with maximal torus \(T=B_+\cap B_-\).
Let \(U\) be the unipotent radical of \(B\), and denote the quotient
maps by \(q_\pm:B_\pm\to T\). Choosing an invariant nondegenerate
bilinear form on \(\g\), an orthonormal basis \((h_i)\) of
\(\Lie(T)\), and dual root vectors \(E_{\pm\alpha}\) for positive
roots \(\alpha\) determined by \(B\), we obtain the standard quasi-triangular $r$-matrix
\[
 r_{\mathrm{st}}=\tfrac12\sum_i h_i\otimes h_i+
             \sum_{\alpha>0}E_{-\alpha}\otimes E_\alpha,\qquad
 \Lambda_{\mathrm{st}}
       =\tfrac12(r_{\mathrm{st}}-r_{\mathrm{st}}^{21}),
\]
where the superscript \(21\) interchanges the tensor factors. We use the Poisson tensor
\(\pi_G=\Lambda_{\mathrm{st}}^R-\Lambda_{\mathrm{st}}^L\), following the convention in
\cite[\S2.1]{Morita}, this amounts to using \(r=-r_{\mathrm{st}}\) in
\eqref{mixed}. The projection \(G\to G/B\) determines a Poisson structure on the
quotient that we denote by \(\pi_{G/B}\). Since the symmetric
part of $r$ is nondegenerate, we obtain that the double is the pair groupoid \(D=G\times G\rightrightarrows G\).
With \(r=-r_{\mathrm{st}}\), the endpoint map \((p_+,p_-)\)
takes the pairing given by evaluation to the one given by the negative of the chosen form on
the first copy of \(\g\) and its positive on the second copy.

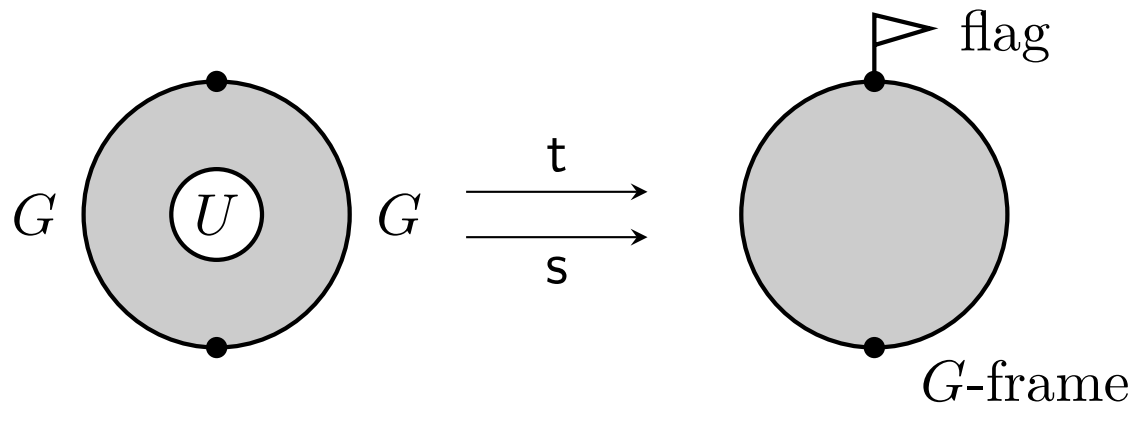
\begin{figure}[htbp]
\centering
\begin{tikzpicture}[x=1cm,y=1cm,>=stealth]
  \draw[thick,fill=black!20] (0,0) circle (.88);
  \draw[thick,fill=white] (0,0) circle (.30);
  \node at (0,0) {\(U\)};
  \fill (0,.88) circle (2pt);
  \fill (0,-.88) circle (2pt);
  \node[left] at (-.92,0) {\(G\)};
  \node[right] at (.92,0) {\(G\)};
  \draw[->] (1.65,.15)--node[above]{\(\tgt\)}(2.85,.15);
  \draw[->] (1.65,-.15)--node[below]{\(\src\)}(2.85,-.15);
  \draw[thick,fill=black!20] (4.35,0) circle (.88);
  \fill (4.35,.88) circle (2pt);
  \fill (4.35,-.88) circle (2pt);
  \draw[thick] (4.35,.88)--(4.35,1.32)
    --(4.72,1.23)--(4.35,1.12);
  \node[right] at (4.78,1.18) {flag};
  \node[right] at (4.52,-1.1) {\(G\)-frame};
\end{tikzpicture}
\caption{The flag annulus and its base. The two maps retain the flag
in the target frame and in the source frame, respectively.}
\label{flag}
\end{figure}

The groupoid integrating the Lie algebroid corresponding to the $\g$-action on $G/B$ is the
decorated annulus of \cite[Ex.~5.26]{PartI},
\cite[proof of Prop.~3.1]{Morita}. Explicitly, the space of arrows is \(\HH_B=G\times(G\times_B U)\), where \(G\times_B U\) is the
quotient by \(b\cdot(g,u)=(gb^{-1},bub^{-1})\). The multiplication of 
two arrows \(h_1=(a_1,[g_1,u_1])\) and \(h_2=(a_2,[g_2,u_2])\) is given by
\[
  (a_1,[g_1,u_1])\,(a_2,[a_1^{-1}g_1,u_2])
  = (a_1a_2,[g_1,u_1u_2]);
\]
choosing representatives with \(g_2=a_1^{-1}g_1\). An element \((a,[g,u])\) has source, target
and moment map
\[
 \src(a,[g,u])=a^{-1}gB,\qquad
 \tgt(a,[g,u])=gB,\qquad
 \Phi_B(a,[g,u])=(gug^{-1}a,a).
\]
 The lagrangian structure on \(\Phi_B\)
is the two-form \(\beta_B\) determined as follows. Let \(\kappa\) denote the chosen invariant bilinear form on \(\g\), with the pairing on \(D=G\times G\) being \((-\kappa,\kappa)\).
Consider the quotient map
\[
 q:G\times G\times U\longrightarrow\HH_B,\qquad
 q(a,g,u)=(a,[g,u]),
\]
and put \(C=gug^{-1}\). Then we may define $\beta_B$ via
\[
 q^*\beta_B
 =
 \frac12\kappa\!\left(
 C^{-1}dC\wedge(dg\,g^{-1}-da\,a^{-1})
 \right)
 -
 \frac12\kappa\!\left(
 g^{-1}dg\wedge du\,u^{-1}
 \right).
\]
Graphically, we may depict $\HH_B$ as the space of representations of the fundamental groupoid of the annulus. One boundary carries a flag with its monodromy lying in the unipotent radical of the stabilizer of that flag while its other boundary has framings on two marked points. In the disc of
Figure~\ref{flag} we have a flag over one marked point and a framing on the other, parallel transport from the flag to the framed point identifies the flag with a point of \(G/B\). Changing that frame
gives the residual \(G\)-action.

The Poisson structure on $G/B$ is integrated by the restriction
\[
 \Phi_B^{-1}(G^*)\rightrightarrows G/B,\qquad
 G^*=\{(b_+,b_-)\in B_+\times B_-:q_+(b_+)q_-(b_-)=1\}.
\]
Here \(G^*\) denotes the
indicated embedded integration of the dual Lie algebra. The residual
moment map in \cite[proof of Prop.~3.1]{Morita} lists the two
holonomies in the reverse order.

\paragraph[Example. Integration of flag configurations]{Integration of flag configurations.}
For \(m\geq1\), we denote by \(\widetilde F_m=G\times_B\cdots\times_B G\) the quotient of \(G^m\) by the intermediate \(B^{m-1}\)-action.
At the \(j\)-th position, we have that \(b\in B\) sends \((g_j,g_{j+1})\)
to \((g_jb,b^{-1}g_{j+1})\). Put \(F_m=\widetilde F_m/B\) for the quotient of the remaining right action on the last factor. Using the expressions \(p_j=g_1\cdots g_j\) for all $j=1\dots m+1$, we obtain the two identifications
\[
 \begin{aligned}
 F_m&\longrightarrow(G/B)^m,&
 [g_1,\ldots,g_m]&\longmapsto(p_1B,\ldots,p_mB),\\
 \widetilde F_{m+1}&\longrightarrow(G/B)^m\times G,&
 [g_1,\ldots,g_{m+1}]&\longmapsto(p_1B,\ldots,p_mB,p_{m+1}).
 \end{aligned}
\]
In these coordinates, the quotient Poisson structures are iterated mixed products \cite[\S2.2 and Thm.~2.1]{Morita},
\cite[Prop.~A.6]{LMY}.

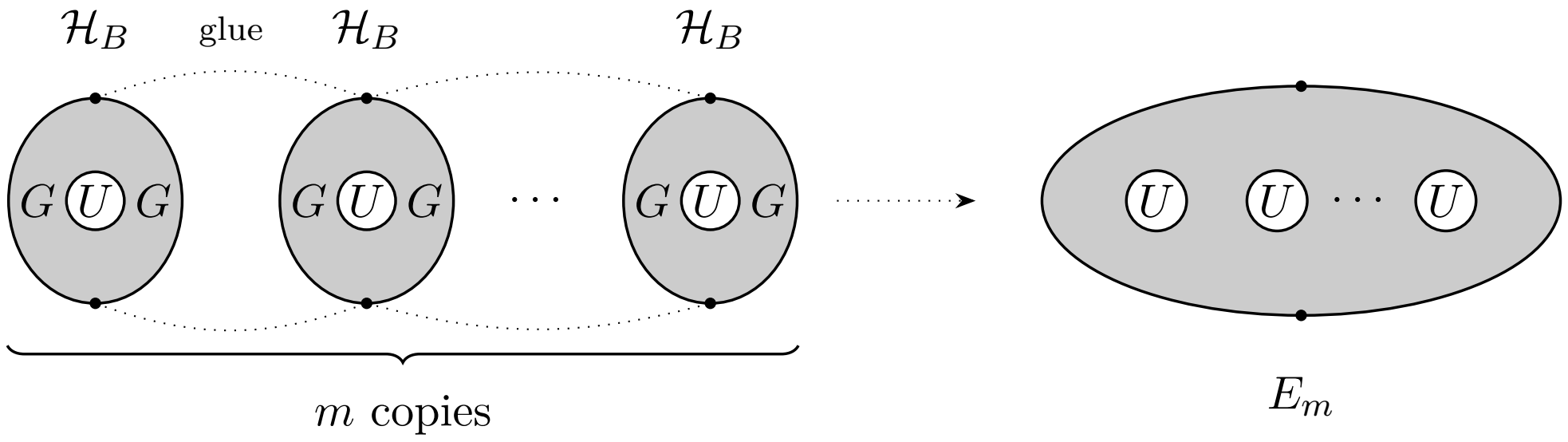
\begin{figure}
    \centering
    \begin{tikzpicture}[line width=.65pt,>=Stealth,font=\normalsize]
  \foreach \x in {0,2.25,5.1} {
    \begin{scope}[shift={(\x,0)}]
      \draw[fill=black!20] (0,0) ellipse (.72 and .85);
      \draw[fill=white] (0,0) circle (.24);
      \node at (0,0) {$U$};
      \node at (-.48,0) {$G$};
      \node at (.48,0) {$G$};
      \fill (0,.85) circle (1.35pt);
      \fill (0,-.85) circle (1.35pt);
      \node at (0,1.42) {$\mathcal H_B$};
    \end{scope}
  }
  \node at (3.68,0) {$\cdots$};
  \draw[dotted,thin] (0,.85) to[out=20,in=160]
    node[above=2pt,font=\scriptsize] {glue} (2.25,.85);
  \draw[dotted,thin] (0,-.85) to[out=-20,in=-160] (2.25,-.85);
  \draw[dotted,thin] (2.25,.85) to[out=15,in=165] (5.1,.85);
  \draw[dotted,thin] (2.25,-.85) to[out=-15,in=-165] (5.1,-.85);
  \draw[decorate,decoration={brace,mirror,amplitude=4pt}]
    (-.73,-1.2) -- (5.83,-1.2)
    node[midway,below=8pt] {$m\text{ copies}$};

  \draw[->,dotted,thin] (6.15,0) -- (7.3,0);

  \draw[fill=black!20] (10,0) ellipse (2.15 and .95);
  \foreach \x in {8.8,9.8,11.2} {
    \draw[fill=white] (\x,0) circle (.25);
    \node at (\x,0) {$U$};
  }
  \node at (10.5,0) {$\cdots$};
  \fill (10,.95) circle (1.35pt);
  \fill (10,-.95) circle (1.35pt);
  \node at (10,-1.62) {$E_m$};
\end{tikzpicture}
    \caption{Successive gluing of decorated annuli integrates the full Lie algebroid of the Poisson action on $F_m$}
    \label{fig:gluing1}
\end{figure}

To integrate \(F_m\), one takes \(m\) copies of $\HH_B$, which is represented by the decorated annulus of
Figure~\ref{flag}. According to Corollary \ref{cor:mix}, one needs to take the fibred product of $m$ such groupoids along matching components of their moment maps. Explicitly, this amounts to the following construction. For each \(h_j\in\HH_B\) we write
\(\Phi_B(h_j)=(a_j,b_j)\in D=G\times G\) and define the fibred product with residual moment map
\[
 E_m=\{(h_1,\ldots,h_m): b_j=a_{j+1}\ (1\leq j<m)\},\qquad
 \Phi_m(h_1,\ldots,h_m)=(a_1,b_m).
\]
Geometrically, these equations describe the gluing of successive annuli along boundary arcs.
The resulting surface is a sphere with \(m+1\) holes. Each of the \(m\)
holes is labelled with a flag and the monodromy around the hole is constrained to lie in the unipotent radical of the stabilizer of the flag. The remaining boundary component carries two marked points and its two arc holonomies serve to define \(\Phi_m\), see Figure~\ref{fig:gluing1}. 

By Theorem~\ref{globalprop} (or Example \ref{ex:poiact}), a symplectic
groupoid integrating the mixed Poisson structure is
\[
             \Phi_m^{-1}(G^*)\ \rightrightarrows\ F_m .
\]
For \(m=1\) this gives the groupoid over \(G/B\) mentioned above.

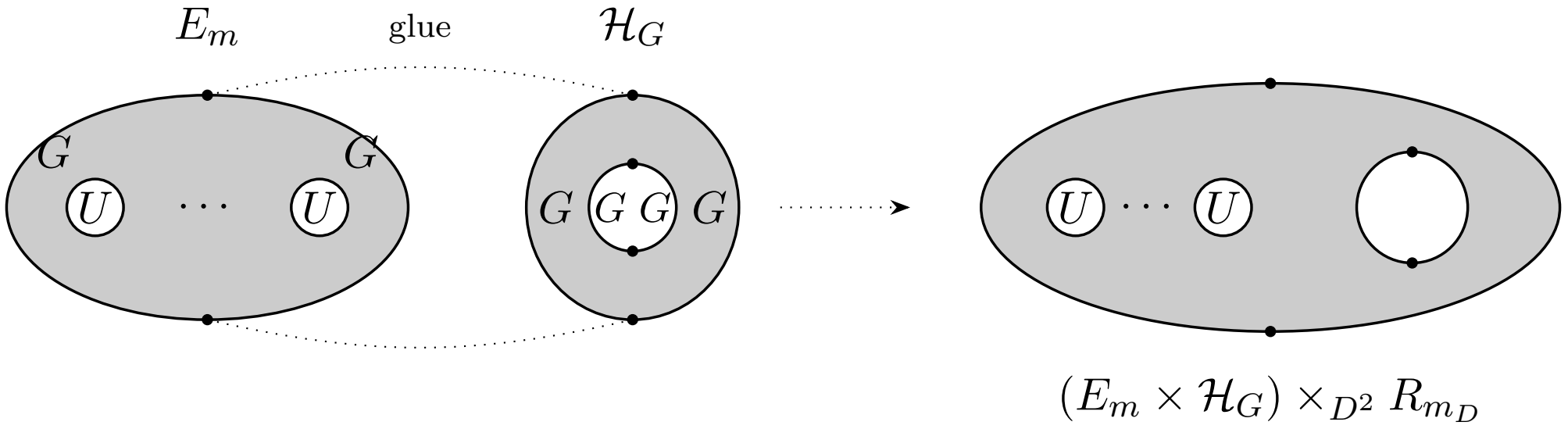
\begin{figure}
    \centering
    \begin{tikzpicture}[line width=.65pt,>=Stealth,font=\normalsize]
  \draw[fill=black!20] (0,0) ellipse (1.7 and .95);
  \foreach \x in {-.95,.95} {
    \draw[fill=white] (\x,0) circle (.24);
    \node at (\x,0) {$U$};
  }
  \node at (0,0) {$\cdots$};
  \fill (0,.95) circle (1.35pt);
  \fill (0,-.95) circle (1.35pt);
  \node at (0,1.53) {$E_m$};
  \node at (-1.3,.47) {$G$};
  \node at (1.3,.47) {$G$};

  \begin{scope}[shift={(3.6,0)}]
    \draw[fill=black!20] (0,0) ellipse (.9 and .95);
    \draw[fill=white] (0,0) circle (.37);
    \fill (0,.95) circle (1.35pt);
    \fill (0,-.95) circle (1.35pt);
    \fill (0,.37) circle (1.35pt);
    \fill (0,-.37) circle (1.35pt);
    \node at (-.65,0) {$G$};
    \node at (.65,0) {$G$};
    \node[font=\small] at (-.19,0) {$G$};
    \node[font=\small] at (.19,0) {$G$};
    \node at (0,1.53) {$\mathcal H_G$};
  \end{scope}
  \draw[dotted,thin] (0,.95) to[out=13,in=167]
    node[above=2pt,font=\scriptsize] {glue} (3.6,.95);
  \draw[dotted,thin] (0,-.95) to[out=-13,in=-167] (3.6,-.95);
  \draw[->,dotted,thin] (4.85,0) -- (5.95,0);

  \draw[fill=black!20] (9,0) ellipse (2.45 and 1.05);
  \foreach \x in {7.35,8.6} {
    \draw[fill=white] (\x,0) circle (.24);
    \node at (\x,0) {$U$};
  }
  \node at (7.98,0) {$\cdots$};
  \draw[fill=white] (10.2,0) circle (.47);
  \fill (9,1.05) circle (1.35pt);
  \fill (9,-1.05) circle (1.35pt);
  \fill (10.2,.47) circle (1.35pt);
  \fill (10.2,-.47) circle (1.35pt);
  \node at (9,-1.62) {$(E_m\times \mathcal H_G)\times_{D^2}{R_{m_D}}$};
\end{tikzpicture}
    \caption{Integration of the Lie algebroid encoding the Poisson action on $\widetilde{F}_{m+1}$}
    \label{fig:gluing2}
\end{figure}

In order to integrate \(\widetilde F_{m+1}\), attach to the same \(m\) decorated annuli the
unrestricted annulus $\HH_G$ integrating the last \(G\)-factor, which is described in
\cite[Prop.~5.23]{PartI}. This gives a sphere with \(m+2\) holes, see Figure~\ref{fig:gluing2}.

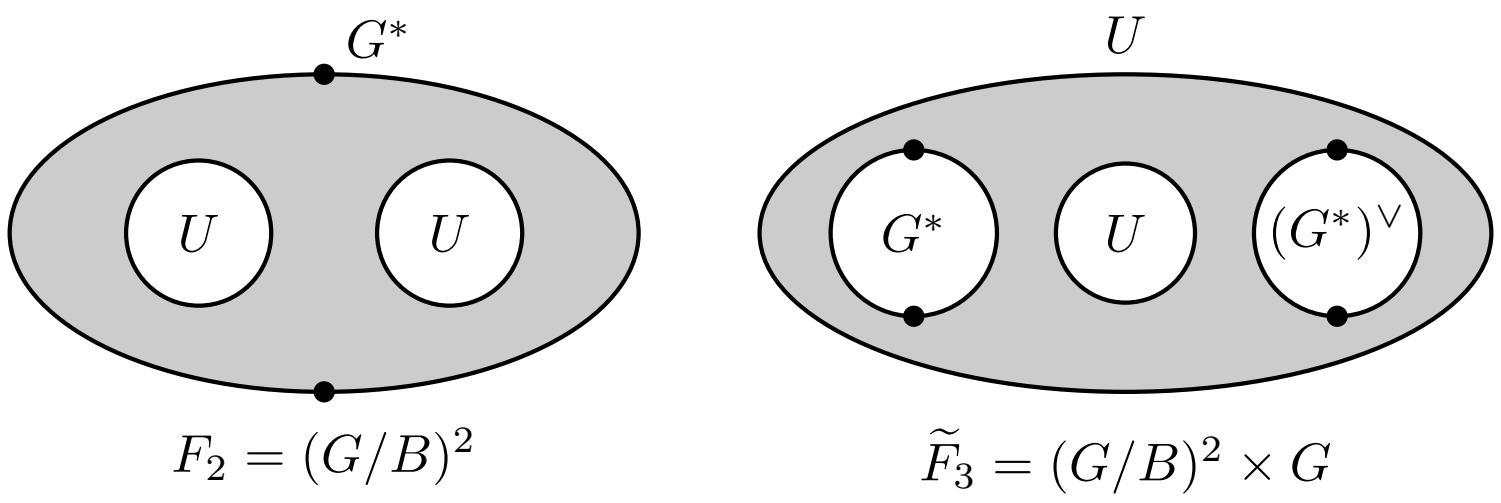
\begin{figure}[htbp]
\centering
\begin{tikzpicture}[x=1cm,y=1cm,thick,font=\small]
  \begin{scope}
    \draw[fill=black!20] (0,0) ellipse (2.08 and 1.05);
    \foreach \x in {-0.83,0.83}{
      \draw[fill=white] (\x,0) circle (.48);
      \node at (\x,0) {\(U\)};
    }
    \foreach \y in {-1.05,1.05}{\fill (0,\y) circle (2pt);}
    \node[above right] at (0,1.02) {\(G^*\)};
    \node[below] at (0,-1.15) {\(F_2=(G/B)^2\)};
  \end{scope}
  \begin{scope}[shift={(5.3,0)}]
    \draw[fill=black!20] (0,0) ellipse (2.42 and 1.05);
    \node[above] at (0,1.05) {\(U\)};
    \draw[fill=white] (0,0) circle (.46);
    \node at (0,0) {\(U\)};
    \foreach \x in {-1.40,1.40}{
      \draw[fill=white] (\x,0) circle (.55);
      \foreach \y in {-.55,.55}{\fill (\x,\y) circle (2pt);}
    }
    \node at (-1.40,0) {\(G^*\)};
    \node at (1.40,0) {\((G^*)^\vee\)};
    \node[below] at (0,-1.15) {\(\widetilde F_3=(G/B)^2\times G\)};
  \end{scope}
\end{tikzpicture}
\caption{Surfaces for the two symplectic groupoids with two flag factors.
Each \(U\) denotes a flag boundary with unipotent monodromy. On each
boundary with two dots the subgroup label constrains the ordered pair
of arc holonomies.}
\label{flag-spheres}
\end{figure}

There are now two free boundaries with marked points at which to impose suitable lagrangian conditions. Write \((G^*)^\vee=\{(v,u):(u,v)\in G^*\}\). Following the convention of \cite[Rem.~3.2]{Morita}, these two ordered
pairs of arc holonomies are constrained to lie in \(G^*\) and \((G^*)^\vee\), respectively. This constraint determines an integration of \(\widetilde F_{m+1}\). Figure~\ref{flag-spheres} shows the two surfaces for \(m=2\).
The construction with the last \(G\)-factor recovers the symplectic
groupoid of \cite[Prop.~3.1]{Morita}.

\phantomsection
\addcontentsline{toc}{section}{\refname}
\bibliographystyle{lagrangian}
\bibliography{lagrangian_refs}

@incollection{LuMomentum,
  author    = {Lu, Jiang-Hua},
  title     = {Momentum mappings and reduction of {Poisson} actions},
  booktitle = {Symplectic Geometry, Groupoids, and Integrable Systems},
  editor    = {Dazord, Pierre and Weinstein, Alan},
  series    = {Mathematical Sciences Research Institute Publications},
  volume    = {20},
  publisher = {Springer},
  address   = {New York},
  year      = {1991},
  pages     = {209--226},
  doi       = {10.1007/978-1-4613-9719-9_15},
}

@article{AMM,
  author  = {Alekseev, Anton and Malkin, Anton and Meinrenken, Eckhard},
  title   = {{Lie} group valued moment maps},
  journal = {J. Differential Geom.},
  volume  = {48},
  number  = {3},
  year    = {1998},
  pages   = {445--495},
  doi     = {10.4310/jdg/1214460860},
  eprint  = {dg-ga/9707021},
}

@article{AKS,
  author  = {Alekseev, Anton and Kosmann-Schwarzbach, Yvette},
  title   = {{Manin} pairs and moment maps},
  journal = {J. Differential Geom.},
  volume  = {56},
  number  = {1},
  year    = {2000},
  pages   = {133--165},
  doi     = {10.4310/jdg/1090347528},
  eprint  = {math/9909176},
}

@article{MouquinLocal,
  author  = {Mouquin, Victor},
  title   = {Local {Poisson} groupoids over mixed product {Poisson} structures and generalised double {Bruhat} cells},
  journal = {J. Symplectic Geom.},
  volume  = {19},
  number  = {4},
  year    = {2021},
  pages   = {993--1045},
  eprint  = {1908.04044},
}

@article{Morita,
  author  = {{\'A}lvarez, Daniel},
  title   = {Morita equivalences, moduli spaces, and flag varieties},
  journal = {Int. Math. Res. Not.},
  volume  = {2026},
  number  = {5},
  year    = {2026},
  eid     = {rnag039},
  doi     = {10.1093/imrn/rnag039},
  url     = {https://doi.org/10.1093/imrn/rnag039},
  eprint  = {2504.09293},
}

@unpublished{PartII,
  author = {{\'A}lvarez, Daniel and Bursztyn, Henrique and Cueca, Miquel},
  title  = {Quasi-hamiltonian groupoids and moduli spaces of flat bundles over surfaces},
}

@misc{PartI,
  author = {{\'A}lvarez, Daniel and Bursztyn, Henrique and Cueca, Miquel},
  title  = {Shifted lagrangian structures in {Poisson} geometry},
  year   = {2026},
  eprint = {2605.29117},
}

@article{CF,
  author  = {Crainic, Marius and Fernandes, Rui Loja},
  title   = {Integrability of {Lie} brackets},
  journal = {Ann. of Math.},
  volume  = {157},
  number  = {2},
  year    = {2003},
  pages   = {575--620},
  doi     = {10.4007/annals.2003.157.575},
  url     = {https://doi.org/10.4007/annals.2003.157.575},
  eprint  = {math/0105033},
}

@article{EV,
  author  = {Etingof, Pavel and Varchenko, Alexander},
  title   = {Solutions of the quantum dynamical {Yang}--{Baxter} equation and dynamical quantum groups},
  journal = {Comm. Math. Phys.},
  volume  = {196},
  year    = {1998},
  pages   = {591--640},
}

@article{FIP,
  author  = {Fernandes, Rui Loja and Iglesias Ponte, David},
  title   = {Integrability of {Poisson}--{Lie} group actions},
  journal = {Lett. Math. Phys.},
  volume  = {90},
  number  = {1--3},
  year    = {2009},
  pages   = {137--159},
  doi     = {10.1007/s11005-009-0348-x},
  url     = {https://doi.org/10.1007/s11005-009-0348-x},
  eprint  = {0902.3558},
}

@article{FOR,
  author  = {Fernandes, Rui Loja and Ortega, Juan-Pablo and Ratiu, Tudor S.},
  title   = {The momentum map in {Poisson} geometry},
  journal = {Amer. J. Math.},
  volume  = {131},
  year    = {2009},
  pages   = {1261--1310},
}

@article{LiParmentier,
  author  = {Li, Luen-Chau and Parmentier, Serge},
  title   = {On dynamical {Poisson} groupoids {I}},
  journal = {Mem. Amer. Math. Soc.},
  volume  = {174},
  number  = {824},
  year    = {2005},
  doi     = {10.1090/memo/0824},
  url     = {https://doi.org/10.1090/memo/0824},
  eprint  = {math/0209212},
}

@incollection{LBS,
  author    = {Li-Bland, David and {\v S}evera, Pavol},
  title     = {Symplectic and {Poisson} geometry of the moduli spaces of flat connections over quilted surfaces},
  booktitle = {Mathematical Aspects of Quantum Field Theories},
  publisher = {Springer},
  year      = {2015},
  pages     = {343--411},
  eprint    = {1304.0737},
  note      = {Example~1.8 and (1.9) in the preprint are Example~8 and (9) in the book},
}

@article{DiracXu,
  author  = {Liu, Zhang-Ju and Xu, Ping},
  title   = {{Dirac} structures and dynamical {$r$}-matrices},
  journal = {Ann. Inst. Fourier},
  volume  = {51},
  year    = {2001},
  pages   = {835--859},
}

@article{Coupled,
  author  = {Liu, Zhang-Ju and Xu, Ping},
  title   = {Dynamical {$r$}-matrices coupled with {Poisson} manifolds},
  journal = {Prog. Theor. Phys. Suppl.},
  volume  = {144},
  year    = {2001},
  pages   = {133--140},
}

@article{LiuXu,
  author  = {Liu, Zhang-Ju and Xu, Ping},
  title   = {The local structure of {Lie} bialgebroids},
  journal = {Lett. Math. Phys.},
  volume  = {61},
  number  = {1},
  year    = {2002},
  pages   = {15--28},
  doi     = {10.1023/A:1020218813233},
  url     = {https://doi.org/10.1023/A:1020218813233},
  eprint  = {math/0210104},
}

@article{Lu,
  author  = {Lu, Jiang-Hua},
  title   = {{Poisson} homogeneous spaces and {Lie} algebroids associated to {Poisson} actions},
  journal = {Duke Math. J.},
  volume  = {86},
  number  = {2},
  year    = {1997},
  pages   = {261--304},
  doi     = {10.1215/S0012-7094-97-08608-7},
  url     = {https://doi.org/10.1215/S0012-7094-97-08608-7},
  eprint  = {q-alg/9503003},
}

@article{LuMouquin,
  author  = {Lu, Jiang-Hua and Mouquin, Victor},
  title   = {Mixed product {Poisson} structures associated to {Poisson} {Lie} groups and {Lie} bialgebras},
  journal = {Int. Math. Res. Not.},
  volume  = {2017},
  number  = {19},
  year    = {2017},
  pages   = {5919--5976},
  doi     = {10.1093/imrn/rnw189},
  url     = {https://doi.org/10.1093/imrn/rnw189},
  eprint  = {1504.06843},
}

@article{LMY,
  author  = {Lu, Jiang-Hua and Mouquin, Victor and Yu, Shizhuo},
  title   = {Configuration {Poisson} groupoids of flags},
  journal = {Int. Math. Res. Not.},
  volume  = {2023},
  number  = {21},
  year    = {2023},
  pages   = {18035--18107},
  doi     = {10.1093/imrn/rnac321},
  url     = {https://doi.org/10.1093/imrn/rnac321},
  eprint  = {2109.03724},
}

@article{MS,
  author  = {Meusburger, Catherine and Sch{\"o}nfeld, Torsten},
  title   = {Gauge fixing and classical dynamical {$r$}-matrices in {$\mathrm{ISO}(2,1)$}-{Chern}--{Simons} theory},
  journal = {Comm. Math. Phys.},
  volume  = {327},
  year    = {2014},
  pages   = {443--479},
  doi     = {10.1007/s00220-014-1938-8},
  url     = {https://doi.org/10.1007/s00220-014-1938-8},
  eprint  = {1203.5609},
}

@article{Mouquin,
  author  = {Mouquin, Victor},
  title   = {The {Fock}--{Rosly} {Poisson} structure as defined by a quasi-triangular {$r$}-matrix},
  journal = {SIGMA},
  volume  = {13},
  year    = {2017},
  pages   = {063, 13~pages},
  doi     = {10.3842/SIGMA.2017.063},
  url     = {https://doi.org/10.3842/SIGMA.2017.063},
  eprint  = {1610.01782},
}

@article{PTVV,
  author  = {Pantev, Tony and To{\"e}n, Bertrand and Vaqui{\'e}, Michel and Vezzosi, Gabriele},
  title   = {Shifted symplectic structures},
  journal = {Publ. Math. Inst. Hautes {\'E}tudes Sci.},
  volume  = {117},
  year    = {2013},
  pages   = {271--328},
  doi     = {10.1007/s10240-013-0054-1},
  url     = {https://doi.org/10.1007/s10240-013-0054-1},
  eprint  = {1111.3209},
}

@article{PingTriangular,
  author  = {Xu, Ping},
  title   = {Triangular dynamical {$r$}-matrices and quantization},
  journal = {Adv. Math.},
  volume  = {166},
  number  = {1},
  year    = {2002},
  pages   = {1--49},
  doi     = {10.1006/aima.2001.2000},
  url     = {https://doi.org/10.1006/aima.2001.2000},
  eprint  = {math/0005006},
}

@article{Xu,
  author  = {Xu, Xiaomeng},
  title   = {Generalized classical dynamical {Yang}--{Baxter} equations and moduli spaces of flat connections on surfaces},
  journal = {Comm. Math. Phys.},
  volume  = {341},
  year    = {2016},
  pages   = {523--542},
  doi     = {10.1007/s00220-015-2533-3},
  url     = {https://doi.org/10.1007/s00220-015-2533-3},
  eprint  = {1407.2071v4},
}
\end{document}